\documentclass[conference]{IEEEtran}

\usepackage{cite}
\usepackage{amsmath,amssymb,amsfonts}
\usepackage{algorithmic}
\usepackage{graphicx}
\usepackage{textcomp}
\usepackage{xcolor}
\usepackage{url}
\usepackage{subcaption}
\usepackage{url}
\usepackage{tabularx}
\usepackage{supertabular}
\usepackage{multirow}
\usepackage{boxedminipage}
\usepackage{eurosym}
\usepackage{enumitem}   
\usepackage[us]{datetime} 
\usepackage[draft=False]{hyperref} 
\hypersetup{
     colorlinks=true,
     linkcolor=blue,
     filecolor=blue,
     citecolor=black,      
     urlcolor=cyan,
     }

\newcommand{\CREspecifications}{\formatdate{12}{07}{2019}}

\def\BibTeX{{\rm B\kern-.05em{\sc i\kern-.025em b}\kern-.08em
    T\kern-.1667em\lower.7ex\hbox{E}\kern-.125emX}}
\begin{document}

\title{Probabilistic forecasting for sizing in the capacity firming framework}

\author{\IEEEauthorblockN{Jonathan~Dumas\IEEEauthorrefmark{1}, Bertrand~Corn\'elusse\IEEEauthorrefmark{1}, Xavier~Fettweis\IEEEauthorrefmark{1}, Antonello~Giannitrapani\IEEEauthorrefmark{2}, Simone~Paoletti\IEEEauthorrefmark{2}, Antonio~Vicino\IEEEauthorrefmark{2}}
	\IEEEauthorblockA{\IEEEauthorrefmark{1} \textit{Departments of Electrical Engineering, Computer Science, and Geography} \\
		\textit{University of Li\`ege, Belgium}\\
		\{jdumas, xavier.fettweis, bertrand.cornelusse\}@uliege.be
	}
	\IEEEauthorblockA{\IEEEauthorrefmark{2} \textit{Dipartimento di Ingegneria dell'Informazione e Scienze Matematiche} \\
		\textit{Universit\`a di Siena, Italy}\\
		\{giannitrapani, paoletti, vicino\}@dii.unisi.it
	}
}

\IEEEoverridecommandlockouts
\IEEEpubid{\makebox[\columnwidth]{978-1-6654-3597-0/21/\$31.00~\copyright2021 IEEE \hfill} \hspace{\columnsep}\makebox[\columnwidth]{ }}

\maketitle

\IEEEpubidadjcol

\begin{abstract}
This paper proposes a strategy to \textit{size} a grid-connected photovoltaic plant coupled with a battery energy storage device within the \textit{capacity firming} specifications of the French Energy Regulatory Commission. In this context, the sizing problem is challenging due to the two-phase engagement control with a day-ahead nomination and an intraday control to minimize deviations from the planning. 
The two-phase engagement control is modeled with deterministic and stochastic approaches. The optimization problems are formulated as mixed-integer quadratic problems, using a Gaussian copula methodology to generate PV scenarios, to approximate the mixed-integer non-linear problem of the capacity firming.
Then, a grid search is conducted to approximate the optimal sizing for a given selling price using both the deterministic and stochastic approaches. The case study is composed of PV production monitored on-site at the Li\`ege University (ULi\`ege), Belgium. 
\end{abstract}


\newcommand{\SOC}{\ensuremath{s}}
\newcommand{\charge}{\ensuremath{P^\text{cha}}}
\newcommand{\BESScharge}{\ensuremath{p^\text{s,cha}}} 
\newcommand{\discharge}{\ensuremath{P^\text{dis}}}
\newcommand{\BESSdischarge}{\ensuremath{p^\text{s,dis}}} 
\newcommand{\maxcharge}{\ensuremath{\overline{S}}}
\newcommand{\mincharge}{\ensuremath{\underline{S}}}
\newcommand{\chargerate}{\ensuremath{\overline{S_P}}}
\newcommand{\dischargerate}{\ensuremath{\underline{S_P}}}
\newcommand{\retentionRate}{\ensuremath{\eta^{\text{retention}}}}
\newcommand{\chargeEfficiency}{\ensuremath{\eta^{\text{cha}}}}
\newcommand{\dischargeEfficiency}{\ensuremath{\eta^{\text{dis}}}}
\newcommand{\initialCharge}{\ensuremath{S^\text{init}}}
\newcommand{\minEndCharge}{\ensuremath{\underline{S}^\text{end}}}
\newcommand{\maxEndCharge}{\ensuremath{\overline{S}^\text{end}}}
\newcommand{\finalCharge}{\ensuremath{S^{\text{end}}}}

\newcommand{\price}[1]{\ensuremath{\pi^{\text{#1}}}}
\newcommand{\gridSalePrice}{\price{}}

\newcommand{\exportrealizations}{\ensuremath{e}^\text{m}}
\newcommand{\realizations}{\ensuremath{p}^\text{m}} 
\newcommand{\Maxproduction}{\ensuremath{P}^+} 
\newcommand{\Minproduction}{\ensuremath{P}^-} 
\newcommand{\maxExportToGridControl}{\ensuremath{E}_{t}^{\text{cap}}}

\newcommand{\engagement}{\ensuremath{p}^{\star}}  
\newcommand{\Maxengagement}{\ensuremath{P}^{\star,+}} 
\newcommand{\Minengagement}{\ensuremath{P}^{\star,-}} 
\newcommand{\export}{\ensuremath{e}}
\newcommand{\production}{\ensuremath{p}} 
\newcommand{\exportInit}{\ensuremath{e}^\text{ini}}
\newcommand{\PositiveEnergyDeviation}{\ensuremath{\Delta e^+}}
\newcommand{\PositiveDeviation}{\ensuremath{\delta p^+}}  
\newcommand{\NegativeEnergyDeviation}{\ensuremath{\Delta e^-}}
\newcommand{\NegativeDeviation}{\ensuremath{\delta p^-}}  
\newcommand{\nonSteerable}{\ensuremath{P^\text{m}}}
\newcommand{\PVpower}{\ensuremath{p^\text{PV,m}}} 
\newcommand{\nonSteerableMax}{\ensuremath{\underline{P}_c}}
\newcommand{\PVcapacity}{\ensuremath{\underline{P}_c}} 
\newcommand{\PVgeneration} {\ensuremath{p^\text{PV}}} 
\newcommand{\productionForecast} {\ensuremath{\widehat{P}}}
\newcommand{\PVforecast} {\ensuremath{\widehat{p}^\text{PV}}} 

\newcommand{\temperature}{\ensuremath{T}}
\newcommand{\irradiance}{\ensuremath{I}}
\newcommand{\clearskyirradiance}{\ensuremath{I}^\text{cs}}
\newcommand{\clearskyirradianceForecast}{\ensuremath{\widehat{I}}^\text{cs}}
\newcommand{\irradianceForecast}{\ensuremath{\widehat{I}}}
\newcommand{\temperatureForecast}{\ensuremath{\widehat{T}}}

\newcommand{\forallt}{\ensuremath{\forall t \in \mathcal{T}}}
\newcommand{\forallw}{\ensuremath{\forall \omega \in \Omega}}

\section{Introduction}
The capacity firming framework is mainly designed for islands or isolated markets. For instance, the French Energy Regulatory Commission (CRE) publishes capacity firming tenders and specifications\footnote{\url{https://www.cre.fr/}.}. The system considered is a grid-connected photovoltaic (PV) plant with a battery energy storage system (BESS) for firming the PV generation. At the tendering stage, offers are selected based on the electricity selling price. Then, the successful tenderer builds its plant and sells the electricity exported to the grid at the contracted selling price, but according to a well-defined daily engagement and penalization scheme. The electricity to be injected in or withdrawn from the grid must be nominated the day-ahead, and engagements must satisfy ramping power constraints.
The remuneration is calculated a \textit{posteriori} by multiplying the realized exports by the contracted selling price minus a penalty. The deviations of the realized exports from the engagements are penalized through a function specified in the tender. A peak option can be activated into the contract leading to a significant selling price increase during a short period of time defined a \textit{priori}.
Therefore the BESS is required to shift the PV generation during peak hours to maximize revenue.

The literature provides several approaches to incorporate energy storage in this framework.
An optimal power management mechanism for a grid-connected PV system with storage is implemented in \cite{riffonneau2011optimal} using Dynamic Programming (DP) and is compared with simple ruled-based management. The BESS is used to perform peak shaving by day-ahead management of the system with a predictive algorithm. The aging process is modeled and taken into account with a cost in the objective function.
A stochastic dual dynamic programming approach is described and implemented  by \cite{pereira1991multi} for scheduling of a 39-reservoir system. This methodology avoids the well-known "curse of dimensionality" of DP by approximating the expected-cost-to-go functions of stochastic dynamic programming by piecewise linear functions obtained from the dual solutions of the optimization problem at each stage.
In the continuity of the stochastic dynamic approach, the sizing and control of an energy storage system to mitigate wind power uncertainty is addressed by \cite{haessig2014dimensionnement}. The wind farm operator is committed on a day-ahead basis to a production engagement.
A data-fitted autoregressive model that captures the day-ahead forecasts correlation is used to quantify the impact of correlation on storage sizing. Then, a parametric variation based on Monte Carlo simulation enables the quantitative assessment of the impact of forecast error correlation on the storage performance. This study demonstrated that discarding the autocorrelation in such an analysis can underestimate the storage capacity by an order of magnitude.
Finally, three distinct optimization strategies, mixed-integer quadratic programming (MIQP), simulation-based genetic algorithm, and expert-based heuristic are empirically compared by \cite{n2019optimal} in the CRE framework with the two-phase engagement control approach. This study demonstrated that future-aware algorithms, using historical solar irradiation data in place of forecasts, achieved approximately savings of 10 \% of the daily profits in comparison with the business-as-usual expert-based strategy approach.

Notwithstanding these studies, the problem of modeling a two-phase engagement control with a stochastic approach in the CRE capacity framework is still an open issue. In addition, the sizing has only been investigated by \cite{haessig2014dimensionnement} in a similar context but with a wind farm. The main thread motivating the contribution of this paper is to extend the work of \cite{dumas2020stochastic} to address the sizing of a grid-connected PV plant and a BESS subject to grid constraints in this context. 
However, the sizing problem is difficult due to this two-phase engagement control with a day-ahead nomination and an intraday control to minimize deviations from the planning. A single sizing optimization problem would result in a non-linear bilevel optimization problem with an upper level, the sizing part, and a lower level, the two-phase engagement control. Thus, a grid search is conducted to approximate the optimal sizing for a given selling price using both the deterministic and stochastic approaches. 
In the two-phase engagement control, the PV uncertainty is taken into account at the planning stage by considering a stochastic approach that uses PV scenarios generated by a Gaussian copula methodology. The planner determines the engagement profile on a day-ahead basis given a selling price profile, the PV scenarios, and the current state of the system, including the battery state of charge. Then, the engagement plan is sent to the grid operator and to the controller that computes every 15 minutes the production, \textit{i.e.}, injection or withdrawal, and set-points, \textit{e.g.}, BESS charge or discharge, until the end of the day. The optimization problems are formulated as Mixed-Integer Quadratic Problems (MIQP) with linear constraints. 

The main novelties of this paper in the CRE capacity framework context are three-fold.
First, a MIQP formulation is proposed to address the planning stage of the two-phase engagement control, that is compatible with a scenario approach, to approximate the Mixed-Integer Non-Linear Programming problem generated by the CRE non-convex penalty function. It is compared to the deterministic formulation using perfect knowledge of the future and PV point forecasts on empirical data from the PV production monitored on-site at the Li\`ege University (ULi\`ege), Belgium.
Second, a transparent and easily reproducible Gaussian copula methodology is implemented to generate PV scenarios based on the parametric Photovoltaics for Utility Scale Applications (PVUSA) model using a regional climate model.
Finally, the sizing of the system is addressed by a grid search to approximate the optimal sizing for a given selling price using both the deterministic and stochastic approaches. 

The remainder of the paper is organized as follows. Section \ref{sec:capacity_firming_process} details the two-phase engagement control. Section \ref{sec:pb_statement} defines the problem statement. Section \ref{sec:forecasting} provides the PVUSA parametric point forecast model, and the Gaussian Copula approach to generate PV scenarios. Section \ref{sec:sizing} investigates the system sizing using both the deterministic and stochastic MIQP approaches. Finally, Section \ref{sec:conclusions} summarizes the main findings and highlights ideas for further work. 

\section{Capacity firming process}\label{sec:capacity_firming_process}

Sections \ref{sec:engagement_process} and \ref{sec:control_process} present the capacity firming process decomposed into a day-ahead engagement and a real-time control steps. Each day is composed of $T$ planning/controlling time periods, with $t$ the planning/controlling time period index, $\mathcal{T}$ the set of planning/controlling time periods, and $\Delta_t$ the planning/controlling time period duration. 

\subsection{Day-ahead engagement process}\label{sec:engagement_process}

The planner computes on a day-ahead basis, before a deadline, a vector of engagements $[\engagement_1, \cdots,\engagement_T]^\intercal$, based on PV generation forecasts. The engagements are accepted by the grid operator if they satisfy the constraints
\begin{subequations}\label{eq:engagement_constraints}
	\begin{align}
	|\engagement_t-\engagement_{t-1}| & \leq  \Delta_{p,t}^{\star} , \ \forallt \setminus \{1\}	\\
	- \engagement_t  & \leq - \Minengagement_t  , \ \forallt\\
	\engagement_t  & \leq  \Maxengagement_t , \ \forallt ,
	\end{align}
\end{subequations}
with $\Delta_{p,t}^{\star}$ being a power ramping constraint, a fraction of the total installed capacity $\PVcapacity$, determined at the tendering stage and imposed by the grid operator. 

\subsection{Real-time control process}\label{sec:control_process}

A receding-horizon controller computes at each period the set-points from $t$ to the end of the day, production that can be injected or withdrew at the grid connection point, PV generation, BESS charge or discharge, based on forecasts of PV generation and the engagements. Only the set-points of the first period are applied to the system.
The remuneration is calculated \textit{ex-post} based on the realized production $\realizations_t$. For a given control period, the net remuneration $r_t$ of the plant is the gross revenue $ \Delta_t \gridSalePrice_t  \realizations_t$ minus a penalty $c(\engagement_t,\realizations_t)$, with $\gridSalePrice_t$ being the contracted selling price
\begin{align}\label{eq:remuneration}	
r_t = \Delta_t \gridSalePrice_t  \realizations_t  - c(\engagement_t,\realizations_t),  \ \forallt.
\end{align}
The penalty function $c$ depends on the specifications of the tender. The penalty occurs when the production is outside a tolerance $[\production_t^-, \production_t^+]$, with $\production_t^- =\engagement_t - \Delta_p$ and $\production_t^+ =\engagement_t + \Delta_p$ the lower and upper bounds, $\Delta_p$ being a fraction of the total installed capacity.

\section{Problem statement}\label{sec:pb_statement}

The problem statement follows the abstract formulation defined by \cite{dumas2020coordination} where a global microgrid control problem can be defined as operating a microgrid safely and in an economically efficient manner. In the capacity firming context, a two-stage approach is considered with a planner and a controller. The CRE non-convex penalty that leads to a non-linear problem is modeled by a quadratic formulation with linear constraints. In this study, the planner and controller optimize economic criteria, which are only related to active power. The ancillary or grid services are not in the scope of the capacity firming specifications. The BESS degradation is not taken into account. The planner and controller horizons are cropped to twenty-four hours. 



\subsection{Stochastic planner quadratic formulation (S)}\label{sec:planner_formulations}

A stochastic planner with a MIQP formulation and linear constraints is implemented using a scenario-based approach. The planner computes on a day-ahead basis the engagement plan $\engagement_t, \ \forallt$ to be sent to the grid. The objective of the stochastic planner to minimize is
\begin{align}\label{eq:S_objective_1}	
J_S &  = \mathop{\mathbb{E}} \bigg[ \sum_{t\in \mathcal{T}} - \Delta_t \gridSalePrice_t \production_t +  c(\engagement_t, \production_t) \bigg] ,
\end{align}
where the expectation is taken with respect to the PV generation $\PVforecast_t$, and $\production_t$ is the power at the coupling point.
Using a scenario-based approach, (\ref{eq:S_objective_1}) is approximated by
\begin{align}\label{eq:S_objective_2}	
J_S &  =  \sum_{\omega \in \Omega} \alpha_\omega  \sum_{t\in \mathcal{T}} \bigg[- \Delta_t \gridSalePrice_t \production_{t,\omega} +  c(\engagement_t, \production_{t,\omega})  \bigg],
\end{align}
with $\alpha_\omega$ the probability of scenario $\omega\in \Omega$, and $\sum_{\omega \in \Omega} \alpha_\omega = 1$. The optimization problem is
\begin{subequations}
\label{eq:S_formulation}	
\begin{align}
\min  & \sum_{\omega \in \Omega} \alpha_\omega  \sum_{t\in \mathcal{T}} \bigg[ - \Delta_t \gridSalePrice_t  \production_{t,\omega} +  c^-_{t, \omega}   \bigg] \quad \text{s.t.:}  \\
c^-_{t, \omega}  & = \frac{\Delta_t \gridSalePrice_t }{\PVcapacity} \NegativeDeviation_{t, \omega}  \big( \NegativeDeviation_{t, \omega}  + 4 \Delta_p\big) , \forallt \\
- \NegativeDeviation_{t, \omega} & \leq   \production_{t, \omega} -\production_t^-   , \forallt \label{eq:S_underprod_cst}	\\
\production_{t, \omega} &  \leq \production_t^+  , \forallt \label{eq:S_overprod_cst}	 \\
\engagement_t-\engagement_{t-1} & \leq  \Delta_{p,t}^{\star},  \forallt \setminus \{1\} \label{eq:S_engagement_cst1}\\
\engagement_{t-1}-\engagement_t & \leq   \Delta_{p,t}^{\star},  \forallt \setminus \{1\} \\
\engagement_t &  \leq \Maxengagement_t,  \forallt \\
- \engagement_t &  \leq -  \Minengagement_t,  \forallt \label{eq:S_engagement_cst2}\\
\PVgeneration_{t,\omega} & \leq \PVforecast_{t,\omega}, \forallt \label{eq:S_action_cst1}		\\
\BESScharge_{t,\omega} & \leq b_{t,\omega} \chargerate , \forallt \\
\BESSdischarge_{t,\omega} &\leq (1-b_{t,\omega})\dischargerate  , \forallt  \\
- \SOC_{t,\omega} & \leq -\mincharge , \forallt\\
\SOC_{t,\omega} & \leq \maxcharge , \forallt \label{eq:S_action_cst2}	\\
\production_{t,\omega} &  = \PVgeneration_{t,\omega} + \BESSdischarge_{t,\omega} - \BESScharge_{t,\omega}  , \forallt \label{eq:S_energy_flows}	\\
\production_{t,\omega} &  \leq \Maxproduction_t , \forallt \label{eq:S_production1}\\
- \production_{t,\omega}&  \leq - \Minproduction_t , \forallt \label{eq:S_production2}\\
\SOC_{1,\omega} - \initialCharge & =  \Delta_t (  \chargeEfficiency  \BESScharge_{1,\omega} - \frac{\BESSdischarge_{1,\omega}}{\dischargeEfficiency}   )  \label{eq:S_soc_dynamic1}\\
\SOC_{t,\omega} - \SOC_{t-1, \omega} & =  \Delta_t (  \chargeEfficiency  \BESScharge_{t,\omega} - \frac{\BESSdischarge_{t,\omega}}{\dischargeEfficiency}   ),  \forallt \setminus \{1\}\\
\SOC_{T,\omega} & = \finalCharge = \initialCharge. \label{eq:S_soc_dynamic2}
\end{align}
\end{subequations}
The optimization variables are $\engagement_t$ (engagement at the coupling point), $\production_{t,\omega}$ (production at the coupling point), $\forallw$, $\NegativeDeviation_{t,\omega}$ (underproduction), $\PVgeneration_{t,\omega}$ (PV generation), $b_{t,\omega}$ (binary variable), $\BESScharge_{t,\omega}$ (BESS charging power), $\BESSdischarge_{t,\omega}$ (BESS discharging power), and $\SOC_{t,\omega}$ (BESS state of charge).
The CRE non-convex piecewise quadratic penalty function is modeled by the constraints (\ref{eq:S_underprod_cst}) $\forallw$, that defines the variables $\NegativeDeviation_{t, \omega}  \in \mathbb{R}_+$ to model the quadratic penalty for underproduction, and (\ref{eq:S_overprod_cst}) $\forallw$ forbidding overproduction that is non-optimal as curtailment is allowed \cite{n2019optimal}.
From (\ref{eq:engagement_constraints}), the engagement constraints are (\ref{eq:S_engagement_cst1})-(\ref{eq:S_engagement_cst2}), where the ramping constraint on $\engagement_1$ is deactivated to decouple consecutive days of simulation. The set of constraints that bound $\PVgeneration_{t,\omega}$, $\BESScharge_{t,\omega}$, $\BESSdischarge_{t,\omega}$, and $\SOC_{t,\omega}$ variables are (\ref{eq:S_action_cst1})-(\ref{eq:S_action_cst2}) $\forallw$ where $\PVforecast_{t,\omega}$ are PV scenarios, and $b_{t,\omega}$ are binary variables that prevent the BESS from charging and discharging simultaneously. The power balance equation and the production constraints are (\ref{eq:S_energy_flows}) and (\ref{eq:S_production1})-(\ref{eq:S_production2}) $\forallw$. 
The dynamics of the BESS state of charge is provided by constraints (\ref{eq:S_soc_dynamic1})-(\ref{eq:S_soc_dynamic2}) $\forallw$. Note, the parameters $\finalCharge$ and $\initialCharge$ are introduced to decouple consecutive days of simulation.

The deterministic (D) formulation of the planner is a specific case of the stochastic formulation by considering only one scenario where $\PVgeneration_{t,\omega}$ become $\PVforecast_t$, PV point forecasts. The deterministic formulation with perfect forecasts ($D^\star$) is D with $\PVforecast_t = \PVpower_t$ $\forallt$. For both the deterministic planners D and $D^\star$, the optimization variables are $\engagement_t$, $\production_t$, $\NegativeDeviation_t$, $\PVgeneration_t$, $b_t$, $\BESScharge_t$, $\BESSdischarge_t$, and $\SOC_t$.

\subsection{Oracle controller}\label{sec:controller_formulations}

The oracle controller is an ideal real-time controller that assumes perfect knowledge of the future by using as inputs the engagement profile to compute the set-points, maximize the revenues and minimize the deviations from the engagements. The oracle controller is $D^\star$ where the engagements are parameters.

\section{Forecasting methodology}\label{sec:forecasting}

\subsection{PV point forecast parametric model}

A PV plant can be modeled using the well-known PVUSA parametric model \cite{dows1995pvusa}, which expresses the instantaneous generated power as a function of irradiance and air temperature
\begin{align}\label{eq:PVUSA}	
\PVgeneration &  = aI + bI^2 + cIT ,
\end{align}
where $\PVgeneration$, $I$ and $T$ are the generated power, irradiance and air temperature, respectively, and $a > 0$, $b < 0$, $c < 0$ are the PVUSA model parameters. These parameters are estimated following the algorithm of \cite{bianchini2013model} that efficiently exploits only the power generation measurements and the theoretical clear-sky irradiance, and is characterized by very low computational effort. The same implementation of the algorithm is used with a sliding window of 12 hours. The parameters reached the steady-state values in 50 days with
\begin{align}\label{eq:theta_uliege_values}	
\theta &  = [0.573, -7.68 \cdot 10^{-5}, -1.86 \cdot 10^{-3}]^T .
\end{align}
%
The weather hindcasts from the MAR regional climate model \cite{fettweis2017reconstructions}, provided by the Laboratory of Climatology of the Li\`ege University, are used as inputs of the parametric PVUSA model. The MAR regional climate is forced by the ERA5 reanalysis database, the fifth-generation European Centre for Medium-Range Weather Forecasts atmospheric reanalysis of the global climate, to produce weather hindcasts. 

\subsection{Gaussian copula-based PV scenarios }

The Gaussian copula approach has been used to generate wind and PV scenarios \cite{pinson2009probabilistic,golestaneh2016generation}. To the best of our knowledge, there is almost no guidance available on which copula family can describe correlated variations in PV power \cite{golestaneh2016generation}, hence the Gaussian copula family is selected instead of copulas like Archimedean or Elliptical.

In this section, let $t$ be the current time index, $k$ be the lead time of the prediction, $Z = \{Z_1, ..., Z_T \}$ a multivariate random variable, $F_{Z_k}(\cdot)$, $k = 1, ..., T $ the marginal cumulative distribution functions, and $R_Z$ the correlation matrix. The goal is to generate samples of $Z$. The Gaussian copula methodology consists of generating a sample $g = \{g_1, ..., g_T \}$ from the Normal distribution $\mathcal{N}(0,R_Z)$. Then, to transform each entry $g_k$ through the standard normal cumulative distribution function $\phi(\cdot)$: $u_k = \phi(g_k)$, $k = 1, ..., T $, and finally, to apply to each entry $u_k$ the inverse marginal cumulative distribution function of $Z_k$: $z_k = F_{Z_k}^{-1}(u_k)$, $k = 1, ..., T $.
In our case, $Z$ is the error between the PV measurement and the point forecast 
\begin{align}\label{eq:Z_definition}
Z_k = \PVpower_k - \PVforecast_k \quad k = 1, ..., T .
\end{align}
$\widehat{F}_{Z_k}(\cdot)$ and $\widehat{R}_Z$ are estimated from the data. Following the methodology described above, a PV scenario $\omega$ is generated based on a sample $\omega$ of PV errors 
\begin{align}\label{eq:PV_scenario_generation}
\PVforecast_{k,\omega} = \PVforecast_k + z_{k,\omega} \quad k = 1, ..., T .
\end{align}
Finally, the set of PV scenarios is generated using the Gaussian copula approach based on the PVUSA point forecasts. 

\section{Sizing study}\label{sec:sizing}

The goal is to determine the optimal sizing of the BESS and PV for a given selling price and its related net revenue over the lifetime project to optimally bid at the tendering stage. 

\subsection{Problem statement and assumptions}\label{sec:sizing_pb_statement}

\begin{figure}[!htb]
	\centering
	\includegraphics[width=\linewidth]{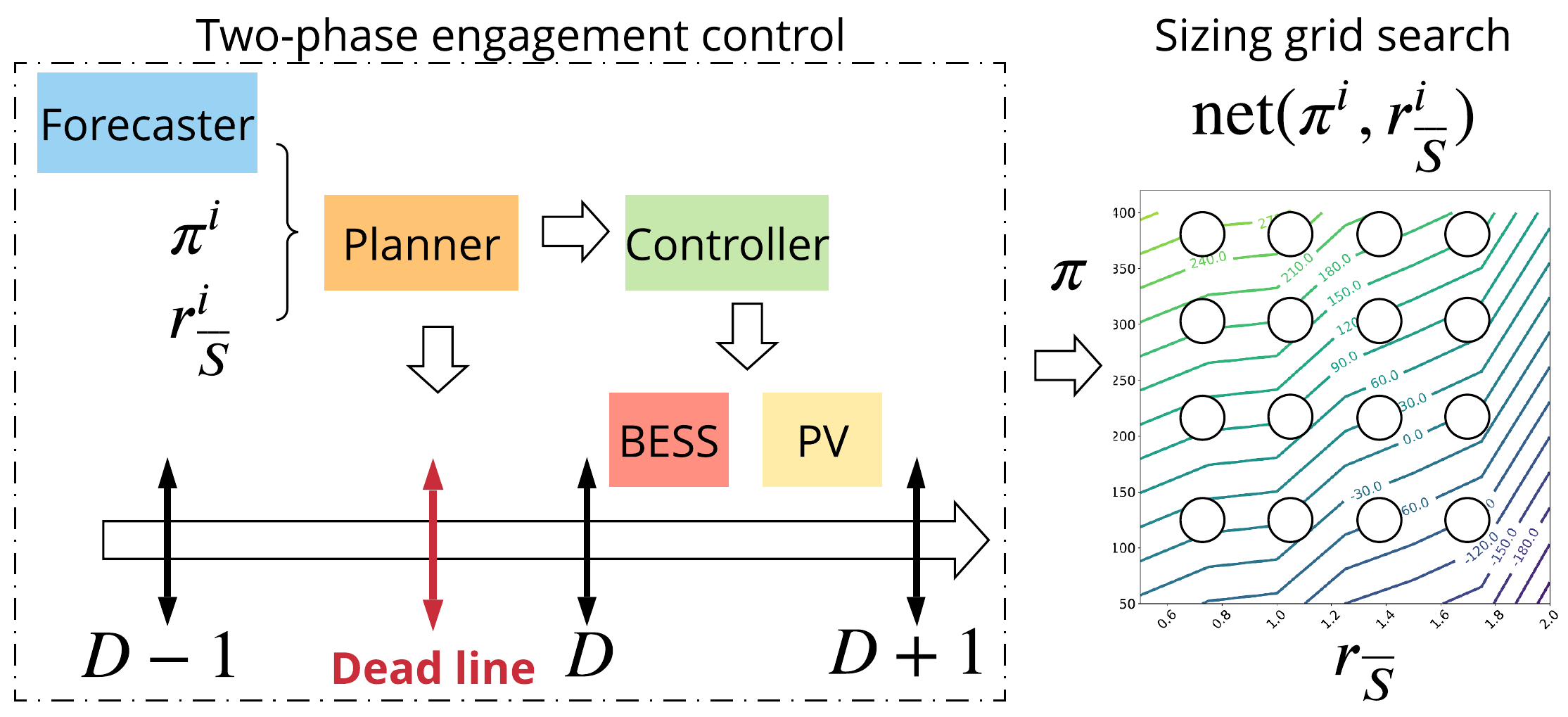}
	\caption{Sizing approach.}
	\label{fig:approach}
\end{figure}

The ratio $r_{\maxcharge} = \frac{\maxcharge}{\PVcapacity}$, BESS maximum capacity over the PV installed capacity, is introduced to model several BESS and PV system configurations. The total exports, imports, deviation costs, number of charging and discharging cycles for several $r_{\maxcharge}$ and selling prices are computed. Based on the BESS and PV CapEx (Capital Expenditure) I and OpEx (Operational Expenditure) O\&M, the LCOE (Levelized Cost Of Energy) is calculated
\begin{align}\label{eq:lcoe_used}
\text{LCOE} & = \frac{\text{CRF}\cdot \text{I} + \text{O\&M} + \text{W} + \text{C}}{\text{E}},
\end{align}
with CRF the Capital Recovery Factor (or Annuity Factor), and E, W, C, the annual export, withdrawal, and deviation costs, respectively. Then, the net revenue over the lifetime project is defined as the annual gross revenue R divided by the annual export to the grid minus the LCOE
\begin{align}\label{eq:net}
\text{net}(\gridSalePrice, r_{\maxcharge}) &: = \frac{R}{E} - \text{LCOE} .
\end{align}
The higher the net revenue, the more profitable the system. Section \ref{sec:LCOE} details the LCOE definition and the assumptions to establish (\ref{eq:lcoe_used}) and (\ref{eq:net}). Finally, the optimal sizing for a given selling price is provided by
\begin{align}\label{eq:optimal}
r_{\maxcharge}^\star (\gridSalePrice)& = \arg \max_{r_{\maxcharge} \in \mathcal{R}_{\maxcharge}} \text{net}(\gridSalePrice, r_{\maxcharge}) ,
\end{align}
with $\mathcal{R}_{\maxcharge}$ the sizing space, and the sizing approach is depicted in Figure (\ref{fig:approach}). 

\subsection{Levelized cost of energy (LCOE)}\label{sec:LCOE}

The LCOE (\euro / MWh) is "the unique break-even cost price where discounted revenues, price times quantities, are equal to the discounted net expenses"\footnote{Annex II Metrics \& Methodology of the AR5 IPCC report.}. 
It means the LCOE is the lower bound of the selling price to be profitable over the lifetime of the project and is defined as
\begin{align}\label{eq:lcoe_def1}
\sum_{t=0}^{n} \frac{\text{E}_t \cdot \text{LCOE}}{(1+i)^t}  & =  \sum_{t=0}^{n} \frac{\text{Expenses}_t}{(1+i)^t},
\end{align}
with $i$ (\%) the discount rate, $n$ (years) the lifetime of the project, and $\text{E}_t$ (MWh) the annual energy at year $t$. Considering energy conversion technologies, the lifetime expenses comprise investment costs I that are the upfront investment (or CapEx) (\euro /kW for PV and \euro /kWh for BESS), operation and maintenance cost (or OpEx) \text{O\&M} (\euro /kW for PV and \euro /kWh for BESS), the annual cost of the energy withdrawn from the grid W (\euro), and the annual deviation cost penalty C (\euro)
\begin{align}\label{eq:lcoe_def2}
\text{LCOE} & =  \frac{\sum_{t=0}^{n} \frac{\text{I}_t + \text{O\&M}_t + \text{W}_t + \text{C}_t}{(1+i)^t}}{\sum_{t=0}^{n} \frac{\text{E}_t}{(1+i)^t}}.
\end{align}
We assume the annual cost of the energy $\text{E}_t$ ($\text{W}_t$) exported (withdrawn), the annual OpEx $\text{O\&M}_t$, and the annual deviation cost $\text{C}_t$ are constant during the lifetime of the project: E, W, C and $\text{O\&M} \ \forall t > 0$. The investment costs I are the sum of all capital expenditures needed to make the investment fully operational discounted to $t=0$. Thus, $\text{I}_t = 0 \ \forall t > 0$ and $\text{I}_0 = \text{I}$. Finally, the system is operational at year $t=1$: $\text{E}_0 = 0$, $\text{O\&M}_0 = 0$, $\text{W}_0 = 0$,  and $\text{C}_0 = 0$. (\ref{eq:lcoe_def2}) becomes
\begin{align}\label{eq:lcoe_def3}
\text{LCOE} & =  \frac{\text{I} + (\text{O\&M} + \text{W} + \text{C}) \sum_{t=1}^{n} \frac{1}{(1+i)^t}}{\text{E}\sum_{t=1}^{n} \frac{1}{(1+i)^t}},
\end{align}
that is re-written
\begin{align}\label{eq:lcoe_def4}
\text{LCOE} & = \frac{\text{CRF}\cdot \text{I} + \text{O\&M} + \text{W} + \text{C}}{\text{E}},
\end{align}
with $\text{CRF} = (\sum_{t=1}^{n} \frac{1}{(1+i)^t})^{-1} = \frac{i}{1-(1+i)^{-n}}$.
The gross revenue is the energy exported to the grid that is remunerated at the selling price. The annual gross revenue $R$ is assumed to be constant over the project lifetime. Then, the LCOE can be compared to R divided by the annual export $E$ to assess the financial viability of the project by calculating (\ref{eq:net}).

\subsection{Case study description}\label{sec:case_study}

The ULi\`ege case study is composed of a PV generation plant with an installed capacity $\PVcapacity$ of 466.4 kW. The simulation dataset is composed of the period from August 2019 to December 2019, 151 days in total. The ULi\`ege PV generation is monitored on a minute basis and is resampled to 15 minutes.

The simulation parameters of the planners and the oracle controller are identical. The planning and controlling periods are $\Delta_t = 15$ minutes. The peak hours are between 7 pm and 9 pm (UTC+0). 
The specifications of the tender AO-CRE-ZNI 2019 published on $\CREspecifications$ define the ramping power constraint on the engagements $\Delta_{p,t}^{\star} = 7.5 \% \PVcapacity$ ($15 \% \PVcapacity$) during off-peak (peak) hours. The lower bound on the engagement is $\Minengagement_t = -5\% \PVcapacity$  ($\Minengagement_t  = 20\% \PVcapacity$) during off-peak (peak) hours. The lower bound on the production is $\Minproduction_t = -5\% \PVcapacity$ ($\Minproduction_t  = 15\% \PVcapacity$) during off-peak (peak) hours. The upper bounds on the engagement and production are $\Maxengagement_t, \Maxproduction_t =  \PVcapacity$. Finally, the engagement deadband is $ \Delta_p = 5\%\PVcapacity $. 

The Python Pyomo\footnote{\url{www.http://www.pyomo.org/}} 5.6.7 library is used to implement the algorithms in Python 3.7. IBM ILOG CPLEX Optimization Studio\footnote{\url{https://www.ibm.com/products/ilog-cplex-optimization-studio}} 12.9 is used to solve all the mixed-integer quadratic optimization problems. Numerical experiments are performed on an Intel Core i7-8700 3.20 GHz based computer with 12 threads and 32 GB of RAM running on Ubuntu\footnote{\url{https://ubuntu.com}} 18.04 LTS. 
The average computation time per optimization problem of the S planner with $|\Omega| = 20$ scenarios is 3 (s) for an optimization problem with 15 000 variables and 22 000 constraints.

\subsection{Sizing parameters}\label{sec:sizing_parameters}

The BESS and PV CAPEX are 300 \euro /kWh and 700 \euro /kW, the BESS and PV OPEX are 1 \% of the CAPEX, the project lifetime is 20 years, and the weighted average cost of capital is 5 \%. The BESS lifetime in terms of full charging and discharging cycles is 3 000.
The BESS is assumed to be capable of fully charging or discharging in one hour $\chargerate = \dischargerate = \maxcharge / 1$, with charging and discharging efficiencies $\chargeEfficiency = \dischargeEfficiency = 0.95$. Each simulation day is independent with a discharged battery at the first and last period to its minimum capacity $\initialCharge = \finalCharge = 10 \% \maxcharge$. The BESS minimum and maximum capacities are 10 \% and 90 \% of $\maxcharge$. 
The sizing space is a grid composed of 56 values with $\mathcal{R}_{\maxcharge} =  \{0.5, 0.75, 1, 1.25, 1.5, 1.75, 2\}$ and $\mathcal{P} =  \{50, 100, 150, 200, 250, 300, 350, 400\}$.

\subsection{Sizing results}\label{sec:sizing_results}

\begin{figure}[!htb]
	\centering
	\begin{subfigure}{.25\textwidth}
		\centering
		\includegraphics[width=\linewidth]{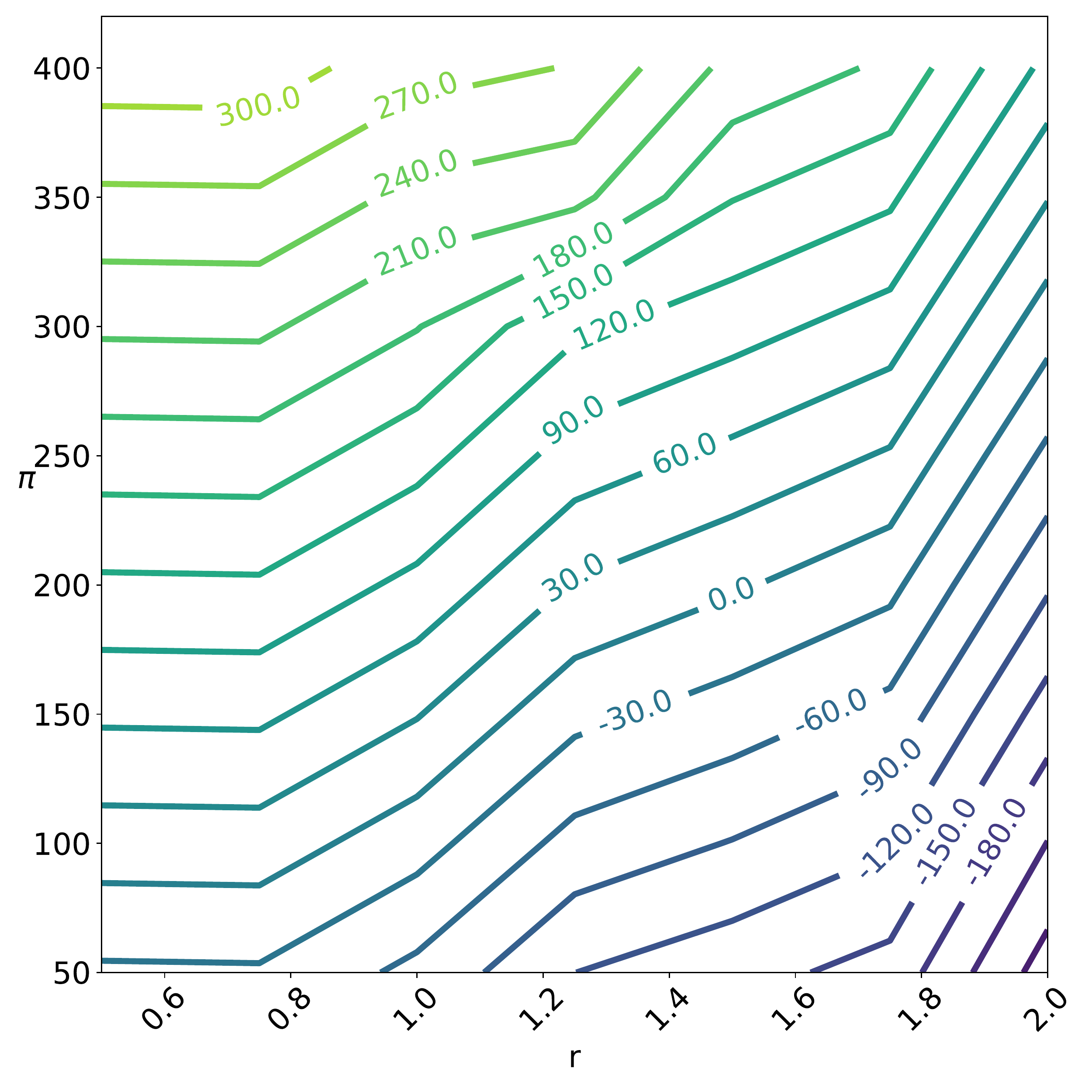}
		\caption{$D^\star$-oracle.}
	\end{subfigure}%
	\begin{subfigure}{.25\textwidth}
		\centering
		\includegraphics[width=\linewidth]{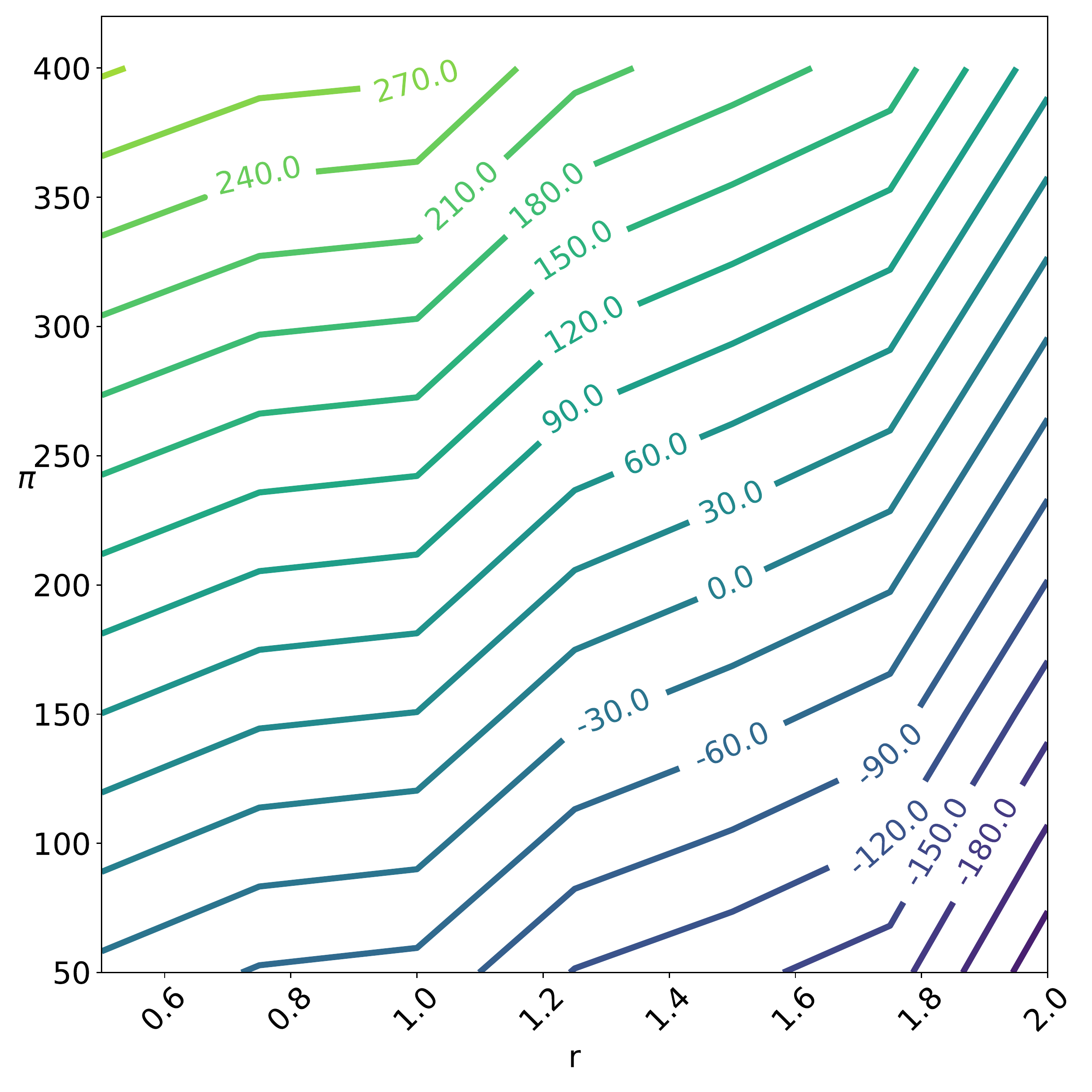}
		\caption{$D$-oracle.}
	\end{subfigure}
	\begin{subfigure}{.25\textwidth}
		\centering
		\includegraphics[width=\linewidth]{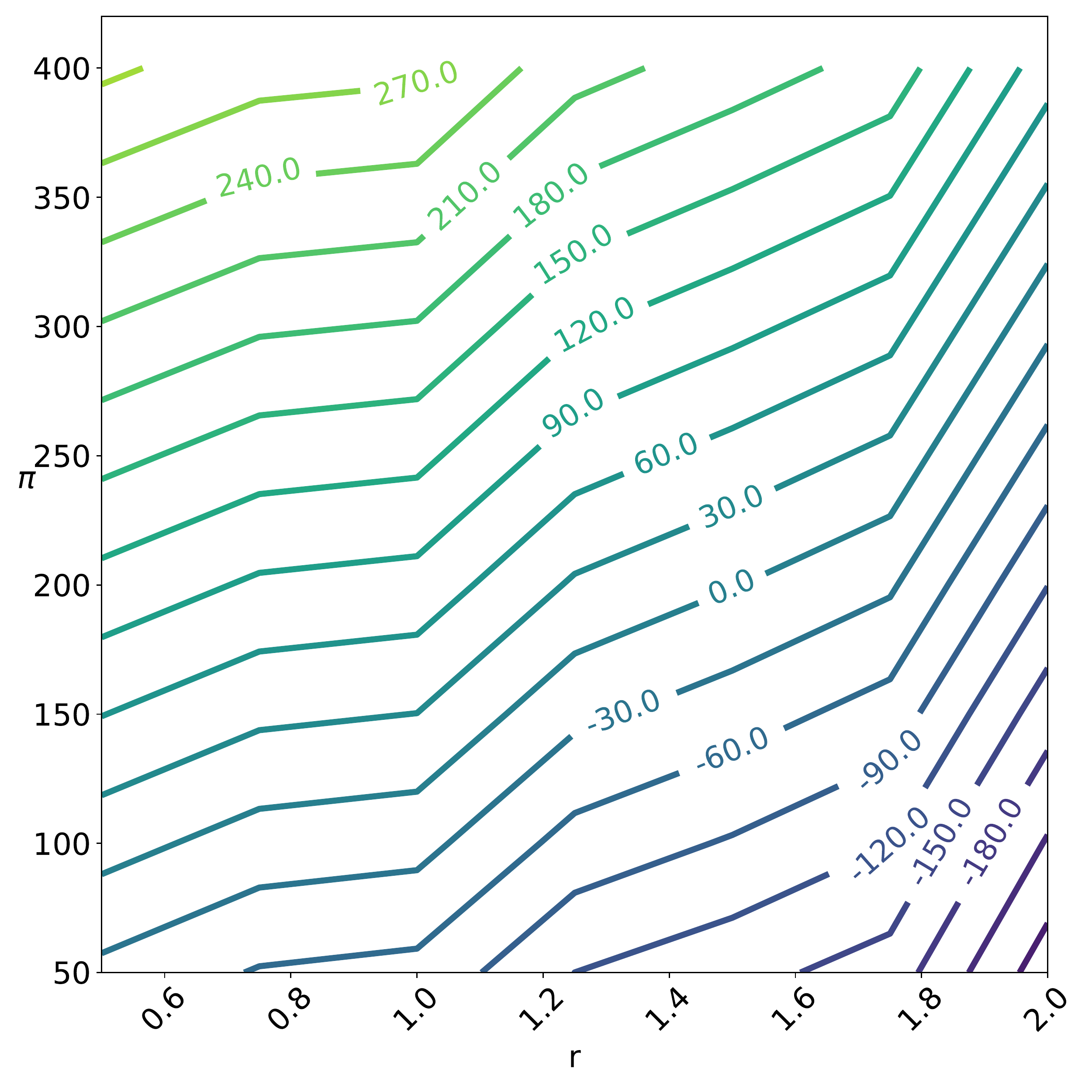}
		\caption{$S^{20}$-oracle.}
	\end{subfigure}%
	\begin{subfigure}{.25\textwidth}
		\centering
		\includegraphics[width=\linewidth]{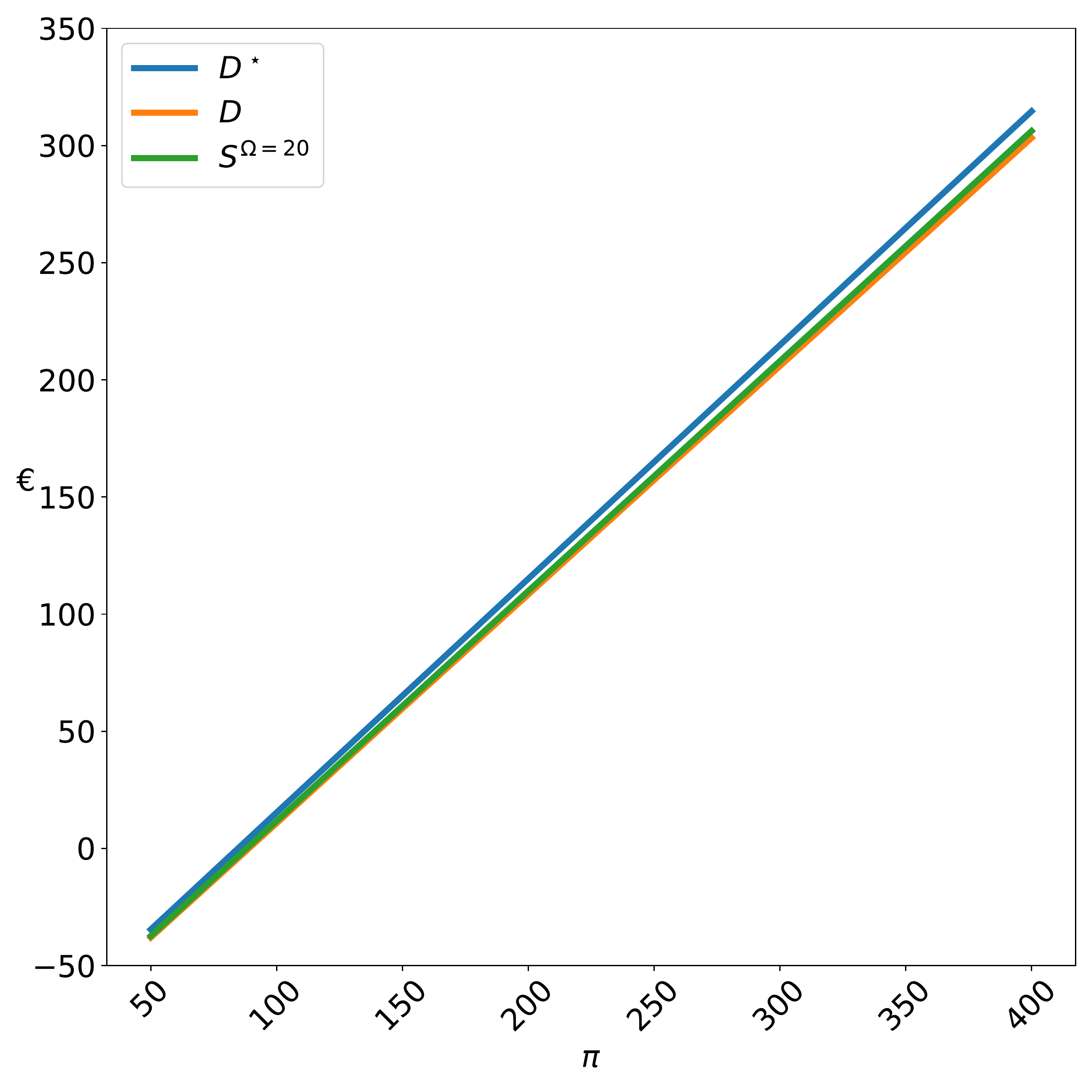}
		\caption{$ \text{net}(\gridSalePrice, r_{\maxcharge}^\star = 0.5)$.}
	\end{subfigure}
	\caption{Sizing results: $\text{net}(\gridSalePrice, r_{\maxcharge})$ (\euro / MWh).}
	\label{fig:sizing_QP_oracle_net}
\end{figure}

The $D^\star$, $D$, and $S^{\Omega = 20}$ planners are used with the oracle controller over the simulation dataset. E, W, C, and the number of full charging and discharging cycles are calculated by extrapolating the results from 151 days to one year.
Figure~\ref{fig:sizing_QP_oracle_net} provides the grid search sizing results for the three planners. For a given selling price, the net is maximal when $r_{\maxcharge} = $ 0.5. When $r_{\maxcharge} $ increases, the BESS is more and more used to withdraw and export during peak hours. It leads to an increase in the number of charging/discharging cycles that implies an increase of the number of BESS required during the project lifetime, and consequently an increase of the BESS CAPEX. As the LCOE is mainly driven by the BESS CAPEX, $\frac{R}{E}$ is not capable of compensating the LCOE increase resulting in a net decrease. The number of charging and discharging cycles is approximately the same for both the $D$ and $S$ planners, independently of the selling price, and rises from 4 700 with $r_{\maxcharge} = $ 0.5 to 13 000 with $r_{\maxcharge} = $ 2. 

The differences between $D^\star$, $D$, and $S$ planners are small. $D^\star$ tends to overestimate the net by underestimating the LCOE (underestimating the deviation, BESS CAPEX, and withdrawal costs) and overestimating $\frac{R}{E}$. 
However, the minimal selling price to be profitable with $r_{\maxcharge} = 0.5$, is approximately 80 \euro / MWh for all planners as shown by Figure~\ref{fig:sizing_QP_oracle_net}. Then, the higher the selling price, the higher the net. In the CRE specifications, the best tender is mainly selected based on the selling price criterion. A trade-off should be reached between the net and the selling price to be selected.

This sizing study seems to indicate that it is not very sensitive to the control policy, \textit{i.e}, deterministic with perfect knowledge, deterministic with point forecasts, and stochastic with scenarios. 
However, it may be dangerous not considering the uncertainty at the sizing stage and could lead to an overestimation of the system performance and an underestimation of the sizing.
Indeed, the net revenues of the two-phase engagement control are similar between planners explaining the small differences in terms of sizing. Two main limitations could explain this result. First, large deviations (15-20 \%) from the engagement plan at the control step occur rarely. Indeed, the oracle controller may compensate for inadequate planning and limits the deviations. A more realistic controller with point forecasts should be considered. Second, when such deviations occur, most of the time they are within the tolerance deadband where there is no penalty. And when the deviations are outside, the penalty is rather small in comparison with the gross revenue. A sensitivity analysis of the numerical settings of the CRE specifications should be performed.

\section{Conclusions and perspectives}\label{sec:conclusions}

The key idea of this paper is to propose a methodology to size the PV and BESS in the context of the capacity firming framework. Indeed, the two-phase engagement control cannot easily be modeled as a single sizing optimization problem. Such an approach would result in a non-linear bilevel optimization problem difficult to solve with an upper level, the sizing part, and a lower level, the two-phase engagement control. 
The two-phase engagement control is decomposed in two steps: computing the day-ahead engagements, then recomputing the set-points in real-time to minimize the deviations from the engagements. The CRE non-convex penalty function is modeled by a threshold-quadratic penalty that is compatible with a scenario approach. The stochastic formulation using a scenario approach is compared to the deterministic formulation. The PV scenarios are generated using a Gaussian copula methodology and PV point forecasts computed with the PVUSA model. 
%
The minimal selling price to be profitable, on this dataset, in the context of the capacity firming framework is approximately 80 \euro / MWh with a BESS having a maximal capacity, fully charged or discharged in one hour, of half the total PV installed power.
The sizing study indicates that it is not very sensitive to the control policy, deterministic with perfect knowledge, deterministic with point forecasts, and stochastic with scenarios. However, further investigations are required to implement a more realistic controller that uses intraday point forecasts, and to conduct a sensitivity analysis on the simulation parameters. 

Several extensions are under investigation. 
A PV generation methodology that is less dependent on the PV point forecasts and takes into account the error dependency of the PV power should be implemented. 
PV scenarios clustering and reduction techniques could be considered to select relevant PV scenarios and improve the stochastic planner results. 
A sizing formulation as a single optimization problem with the PV and BESS capacities as variables to directly compute the optimum and avoid making a grid search. Such a formulation is not trivial due to the specific two-phase engagement control of the capacity firming framework.
Finally, the BESS aging process could be modeled in the sizing study and a dataset with at least a full year of data should be considered to fully take into account the PV seasonality.

\bibliographystyle{ieeetr}
\bibliography{biblio}

\end{document}